\documentclass[10pt,aps,pre,twocolumn,showpacs,superscriptaddress,floatfix,longbibliography]{revtex4-2}

\usepackage{amsmath,amssymb,amsfonts,mathtools}
\usepackage{graphicx}
\usepackage{multirow}
\usepackage{appendix}
\usepackage{xcolor}
\usepackage{siunitx}
\usepackage{yfonts}
\usepackage{esint}
\usepackage{placeins}

\usepackage{url}
\usepackage{hyperref}

\begin{document}

\title{A conservative finite-volume Buckley--Leverett solver with bounded-interval multiwavelet state analysis}

\author{Christian Tantardini}
\email{christiantantardini@ymail.com}
\affiliation{Center for Integrative Petroleum Research, King Fahd University of Petroleum and Minerals, Dhahran 31261, Saudi Arabia.}

\date{\today}

\begin{abstract}
We develop a conservative finite-volume Buckley--Leverett solver equipped with a bounded-interval multiwavelet state-analysis layer. Because non-capillary Buckley--Leverett transport is a nonlinear hyperbolic conservation law with entropy-admissible shocks, the saturation equation is advanced by a conservative finite-volume method with monotone numerical fluxes. The accepted finite-volume state is then embedded in a bounded-domain multiwavelet hierarchy, reconstructed back to cell averages, and used for multiresolution diagnostics. The formal transport accuracy is therefore governed by the underlying finite-volume discretization, while the multiwavelet layer is used for representation, compression, and front-localization diagnostics. Its purpose is instead to quantify whether the deterministic physical-space saturation state can be represented faithfully, compressed in a controlled manner, and used to identify dynamically active front regions. Validation against reference Buckley--Leverett profiles for a Berea benchmark shows accurate saturation histories, spatial profiles, front-position diagnostics, and mass balance. The multiwavelet reconstruction tracks the internal finite-volume state with essentially exact fidelity. Additional thresholding tests show that a substantial fraction of detail coefficients can be discarded while maintaining small reconstruction errors and negligible global mass defect, and fine-level detail activity localizes the moving displacement front. The resulting formulation provides a conservative and reproducible first stage toward future transport-active adaptive multiwavelet solvers for porous-media flow.
\end{abstract}

\maketitle

\section{Introduction}
\label{sec:intro}

Accurate prediction of immiscible two-phase displacement, including waterflooding, water-alternating-gas injection, and polymer-assisted recovery, is central to reservoir engineering because front propagation, breakthrough time, and recovery efficiency directly affect injection and production decisions.\cite{BuckleyLeverett1942,BuckleyLeverett1952,belazreg2019novel}
At the continuum scale, the wetting-phase saturation in incompressible one-dimensional displacement is governed by the Buckley--Leverett equation, obtained from Darcy's law and phase mass conservation.\cite{BuckleyLeverett1942,BuckleyLeverett1952}
In the absence of capillarity, this equation is a nonlinear hyperbolic conservation law and its physically relevant weak solution may contain shocks and rarefactions.
A numerical method must therefore preserve conservative flux balance and propagate fronts with the correct entropy-admissible speed.\cite{kaasschieter1999buckley,leveque2002}
The present work focuses on this non-capillary hyperbolic regime because it isolates the conservative transport mechanism that any multiresolution formulation for Buckley--Leverett flow must respect.

Conservative finite-volume (FV) methods remain a natural baseline for
Buckley--Leverett transport because they update cell averages through
intercell flux differences and therefore encode local conservation
directly.\cite{leveque2002}
Several higher-resolution or adaptive alternatives are available.
Reconstructed FV schemes improve interface states while retaining the
conservative update structure,\cite{leveque2002}
weighted essentially non-oscillatory (WENO) schemes provide high-order
non-oscillatory resolution of steep gradients,\cite{Shu2016high,WENO_methods}
Runge--Kutta discontinuous Galerkin (RKDG) methods use high-order local
polynomial representations for hyperbolic conservation laws,\cite{cockburn1998rkdg,cockburn2001rkdg}
and adaptive mesh-refinement strategies concentrate grid resolution near shocks
and other localized features.\cite{berger1984amr,berger1989amr}
In all cases, the conservative flux structure remains central for
shock-compatible Buckley--Leverett transport.
Multiresolution methods provide a complementary hierarchy for such discretizations.
Harten's multiresolution framework represents conservation-law solutions through nested dyadic grid averages and detail information, allowing small details to be associated with smooth regions and large details with localized activity such as shocks or steep gradients.\cite{harten1995multiresolution,harten1993discrete}
This idea underlies adaptive multiresolution FV schemes, where wavelet-like details are used for compression and refinement while the FV structure controls the conservative update.\cite{cohen2003fully,muller2007fully}
A related use of hierarchy appears in multiwavelet-based DG indicators, where high-level difference coefficients identify troubled cells or discontinuities without applying limiting everywhere.\cite{vuik2014multiwavelet,gerhard2015highorder,huang2020adaptive}

In reservoir-flow applications, multiwavelets have also been used in stochastic
or uncertainty spaces for Buckley--Leverett-type problems, while the
physical-space transport operator remains finite-volume.\cite{pettersson2016stochastic}
The deterministic physical-space role considered here is different.
Here, the accepted conservative FV saturation state is projected into a
bounded-interval multiwavelet hierarchy and reconstructed back to cell averages,
so that the representation can be tested directly in the same variables used by
the transport solver.
This construction separates the shock-compatible transport update from the
multiresolution state analysis: the finite-volume residual controls the
conservative propagation of the front, while the multiwavelet layer is used to
quantify reconstruction fidelity, conservation compatibility, compressibility
under detail-coefficient thresholding, and localization of the moving front by
fine-level activity.
The method is validated against reference Buckley--Leverett profiles generated by the standard fractional-flow construction implemented in the \texttt{pywaterflood} package.\cite{Male2024pywaterflood}

The paper is organized as follows.
Section~\ref{sec:phys-model} introduces the one-dimensional Buckley--Leverett model and the conservative FV formulation.
Section~\ref{sec:optionA} describes the bounded-interval multiwavelet representation of the accepted conservative state, including projection, reconstruction, compression, conservation diagnostics, and adaptive-indicator construction.
Section~\ref{sec:numerics} validates the method against reference Buckley--Leverett solutions and analyzes the multiresolution diagnostics.
Section~\ref{sec:conclusions} summarizes the main results and discusses how the
present representation layer can support later transport-active adaptive
multiwavelet formulations.

\section{One-dimensional Buckley--Leverett model and bounded-domain conservative formulation}
\label{sec:phys-model}

We consider immiscible, incompressible displacement of oil by water in a one-dimensional core of length $L$, with spatial coordinate
\begin{align}
x\in[0,L], \qquad t>0.
\end{align}
Under the classical Buckley--Leverett assumptions of constant porosity $\phi$, negligible capillarity, and constant total Darcy velocity $v$, the wetting-phase saturation $S_w(x,t)$ satisfies the conservation law
\begin{align}
\phi\,\partial_t S_w + \partial_x\!\big(v\,f_w(S_w)\big)=0.
\label{eq:bl-1d-phys}
\end{align}
The initial condition is taken as a uniform state,
\begin{align}
S_w(x,0)=S_{w,\mathrm{init}},
\label{eq:bl-init-phys}
\end{align}
while the injection boundary condition at the inlet is
\begin{align}
S_w(0,t)=S_{w,\mathrm{inj}}.
\label{eq:bl-inlet-phys}
\end{align}
At the outlet, the present formulation employs an outflow closure through the numerical interface flux.

Equation \eqref{eq:bl-1d-phys} is the standard one-dimensional Buckley--Leverett transport model. In the non-capillary regime it is a scalar nonlinear hyperbolic conservation law, so its physically relevant weak solution generally contains shocks and rarefactions and must be interpreted in the entropy sense. For this reason, the numerical formulation must preserve conservative flux balance and propagate discontinuities with the correct speed. These requirements are standard in the theory and numerics of hyperbolic conservation laws and are especially important for Buckley--Leverett-type fronts.\cite{leveque2002}

The nonlinear flux in \eqref{eq:bl-1d-phys} is determined by the water fractional-flow function. In the present work, the phase mobilities are modeled by Corey-type relative permeabilities. Introducing the effective saturation
\begin{align}
S_{we}=\frac{S_w-S_{wc}}{S_{w,\mathrm{inj}}-S_{wc}},
\label{eq:effective-sat-phys}
\end{align}
for the present Buckley--Leverett setting, the injected water saturation is taken as the upper mobile endpoint, so that $S_{w,\mathrm{inj}}=1-S_{or}$ in the default Berea benchmark.
The endpoint-scaled relative permeabilities are written as
\begin{align}
k_{rw}(S_w) &= k_{rw}^0\,S_{we}^{\,n_w},\\
k_{ro}(S_w) &= k_{ro}^0\,(1-S_{we})^{\,n_o},
\label{eq:corey-phys}
\end{align}
and the corresponding phase mobilities are
\begin{align}
\lambda_w(S_w)=\frac{k_{rw}(S_w)}{\mu_w},
\qquad
\lambda_o(S_w)=\frac{k_{ro}(S_w)}{\mu_o},
\label{eq:mobilities-phys}
\end{align}
with $\mu_w$ and $\mu_o$ the water and oil viscosities. Corey-type parametrizations remain among the most widely used analytic constitutive models for two-phase relative permeability and provide a convenient and interpretable setting for one-dimensional Buckley--Leverett validation. \cite{Corey1954,BrooksCorey1964}

The water fractional flow is then
\begin{align}
f_w(S_w)=\frac{\lambda_w(S_w)}{\lambda_w(S_w)+\lambda_o(S_w)}.
\label{eq:fractional-flow-phys}
\end{align}
It is convenient to define the conservative transport flux
\begin{align}
F(S_w)=\frac{v}{\phi}\,f_w(S_w),
\label{eq:flux-phys}
\end{align}
so that \eqref{eq:bl-1d-phys} becomes
\begin{align}
\partial_t S_w + \partial_x F(S_w)=0.
\label{eq:cons-law-phys}
\end{align}

We therefore discretize the bounded interval $[0,L]$ into $N$ control volumes
\begin{align}
I_j=[x_{j-\frac12},x_{j+\frac12}],
\qquad
\Delta x=\frac{L}{N},
\qquad
j=1,\dots,N,
\label{eq:grid-phys}
\end{align}
with cell centers
\begin{align}
x_j=\frac12(x_{j-\frac12}+x_{j+\frac12}).
\end{align}
The numerical unknowns are the cell averages
\begin{align}
\bar S_j(t)=\frac{1}{\Delta x}\int_{I_j} S_w(x,t)\,dx.
\label{eq:cellavg-phys}
\end{align}
Integrating \eqref{eq:cons-law-phys} over each control volume gives the semidiscrete conservative update
\begin{align}
\frac{d}{dt}\bar S_j
=
-\frac{1}{\Delta x}
\left(
\widehat F_{j+\frac12}
-
\widehat F_{j-\frac12}
\right),
\label{eq:semidiscrete-fv-phys}
\end{align}
where $\widehat F_{j+\frac12}$ is a consistent numerical flux evaluated at the interface $x_{j+\frac12}$.
This is the standard finite-volume form for a scalar conservation law and directly encodes local conservation through flux differences. 

Two monotone scalar fluxes are considered in the implementation.
The default choice is the Godunov flux, \cite{Godunov1959,leveque2002}
\begin{align}
\widehat F^{G}(a,b)
=
\begin{cases}
\displaystyle \min_{s\in[a,b]} F(s), & a\le b,\\[1ex]
\displaystyle \max_{s\in[b,a]} F(s), & a>b,
\end{cases}
\label{eq:godunov-phys}
\end{align}
which is naturally tied to the entropy solution of the scalar conservation law and is therefore a particularly appropriate choice for Buckley--Leverett transport. 
A Rusanov flux is also available, \cite{Rusanov1961,leveque2002}
\begin{align}
\widehat F^{R}(a,b)
=
\frac12\big(F(a)+F(b)\big)
-
\frac12\alpha(a,b)\,(b-a),
\label{eq:rusanov-phys}
\end{align}
where $\alpha(a,b)$ is a local wave-speed bound, typically obtained from the derivative of $F$.
The Rusanov flux is more diffusive than Godunov but provides a robust fallback and remains monotone for the scalar problem. 

The inlet condition \eqref{eq:bl-inlet-phys} is incorporated through the left interface flux
\begin{align}
\widehat F_{\frac12}=\widehat F(S_{w,\mathrm{inj}},\bar S_1),
\label{eq:left-flux-phys}
\end{align}
while the outlet is treated by an outflow-type closure
\begin{align}
\widehat F_{N+\frac12}=\widehat F(\bar S_N,\bar S_N).
\label{eq:right-flux-phys}
\end{align}
With these definitions, the semidiscrete system \eqref{eq:semidiscrete-fv-phys} is advanced in time by the second-order strong-stability-preserving Runge--Kutta method, \cite{gottlieb2001ssp}
\begin{align}
\bar{\boldsymbol S}^{(1)}
&=
\bar{\boldsymbol S}^{\,n}
+
\Delta t\,\boldsymbol R(\bar{\boldsymbol S}^{\,n}),
\label{eq:ssprk-stage1-phys}
\\
\bar{\boldsymbol S}^{\,n+1}
&=
\frac12\bar{\boldsymbol S}^{\,n}
+
\frac12\Big(
\bar{\boldsymbol S}^{(1)}
+
\Delta t\,\boldsymbol R(\bar{\boldsymbol S}^{(1)})
\Big),
\label{eq:ssprk-stage2-phys}
\end{align}
where $\boldsymbol R$ denotes the conservative residual from \eqref{eq:semidiscrete-fv-phys}.
The time step is restricted by a CFL condition based on the maximum characteristic speed of the fractional-flow flux, which is standard for explicit finite-volume discretizations of nonlinear hyperbolic equations. 

This finite-volume formulation supplies the conservative transport state that is analyzed in the next section.
The subsequent bounded-interval multiwavelet step is applied to accepted finite-volume states and therefore does not modify the interface flux balance defined above. 

\section{Bounded-interval multiwavelet representation of the conservative Buckley--Leverett state}
\label{sec:optionA}

The finite-volume formulation in Section~\ref{sec:phys-model} provides the conservative transport state.
The multiwavelet component introduced here acts on this accepted state rather than on the flux residual.
Thus the numerical evolution remains controlled by the finite-volume update, while the bounded-interval multiwavelet representation supplies a hierarchical description of the transported saturation profile.

We map the physical interval $x\in[0,L]$ to the normalized coordinate
\begin{align}
\xi=\frac{x}{L}\in[0,1].
\label{eq:normalized-coordinate-optionA}
\end{align}
On this bounded interval we introduce a multiwavelet basis $\{\psi_k(\xi)\}$.
In the present implementation, the bounded-domain multiresolution operations are performed with \texttt{VAMPyR} (Very Accurate Multiresolution Python Routines), using a one-dimensional interval representation on $[0,1]$, order $8$, and projection precision $10^{-7}$.\cite{bast_2023_7967323,battistella_2023_10290360,vampyr_github,vampyr_2024}
The use of a bounded-interval representation is important here because the computational domain is a finite core rather than a periodic or infinite interval.

Let
\begin{align}
\bar{\boldsymbol S}
=
(\bar S_1,\dots,\bar S_N)^{\mathsf T}
\end{align}
denote the cell-average vector after an accepted finite-volume transport step.
The first step is to interpret this vector as a piecewise-constant function on the normalized interval:
\begin{align}
S_h(\xi)=\bar S_j,
\qquad
\xi\in\left[\frac{j-1}{N},\frac{j}{N}\right],
\qquad
j=1,\dots,N.
\label{eq:piecewise-constant-optionA}
\end{align}
This identifies the conservative finite-volume state with a bounded-domain function while preserving its cell-average meaning.

The piecewise-constant state is then projected into the bounded-interval multiwavelet space,
\begin{align}
S_h(\xi)\approx S_h^{\mathrm{MW}}(\xi)
=
\sum_k s_k\,\psi_k(\xi),
\label{eq:mw-state-optionA}
\end{align}
where $s_k$ are the multiwavelet coefficients returned by the bounded-domain projection.
The projected field is then mapped back to the finite-volume variable space by cellwise quadrature,
\begin{align}
\bar S_j^{\,\mathrm{MW}}
=
N\int_{(j-1)/N}^{j/N}
S_h^{\mathrm{MW}}(\xi)\,d\xi,
\qquad
j=1,\dots,N.
\label{eq:mw-cellavg-optionA}
\end{align}
This gives the reconstructed cell-average vector
\begin{align}
\bar{\boldsymbol S}^{\,\mathrm{MW}}
=
(\bar S_1^{\,\mathrm{MW}},\dots,\bar S_N^{\,\mathrm{MW}})^{\mathsf T},
\end{align}
which can be compared directly with the original finite-volume state.

The conservation meaning of this projection--reconstruction step is as follows.
The finite-volume update is locally conservative before projection because it is written as a difference of numerical interface fluxes.
The multiwavelet reconstruction is applied only after an accepted finite-volume step and is not used to compute the intercell fluxes in the reported simulations.
Therefore, the local conservation property of the transport update is unchanged.
When the reconstructed state is used for diagnostics or visualization, its agreement with the finite-volume cell averages is controlled by the projection and quadrature errors.
If a thresholded or otherwise modified reconstructed state were later used as input to the next transport step, exact cell-average preservation or a conservative correction would be required.
This distinction is important: in the calculations reported here, the finite-volume state remains the transported state.

The same cell-average vector also provides a simple dyadic multiresolution diagnostic.
Starting from $\bar{\boldsymbol S}$, repeated coarse--fine splitting gives coarse averages and detail coefficients.
At one dyadic split, two neighboring values are replaced by their average and half-difference; repeating this operation over levels gives detail coefficients $d_{\ell,k}$, where $\ell$ denotes the dyadic level and $k$ the location.
The detail energy at level $\ell$ is defined as
\begin{align}
E_\ell
=
\sum_k d_{\ell,k}^2.
\label{eq:detail-energy-optionA}
\end{align}
The sequence $\{E_\ell\}$ measures how strongly the saturation profile occupies the different dyadic levels.
Smooth regions produce small fine-level details, whereas a steep Buckley--Leverett front generates localized fine-scale activity.
This diagnostic role is consistent with multiresolution finite-volume methods, where details are used for compression and refinement, and with multiwavelet DG indicators, where high-level coefficients are used to detect troubled or discontinuous regions.\cite{harten1995multiresolution,cohen2003fully,vuik2014multiwavelet}

To test the compressibility of the representation, we also perform a detail-thresholding experiment on accepted finite-volume states.
The cell-average vector is decomposed into a coarsest mean coefficient and dyadic details,
\begin{align}
\bar{\boldsymbol S}
\longleftrightarrow
\left\{
\bar S_0,\ d_{\ell,k}
\right\},
\end{align}
where $\bar S_0$ denotes the global mean content.
For a prescribed threshold $\varepsilon_{\rm MR}$, the details are replaced by
\begin{align}
d_{\ell,k}^{\,\varepsilon}
=
\begin{cases}
d_{\ell,k}, & |d_{\ell,k}|\ge \varepsilon_{\rm MR},\\
0, & |d_{\ell,k}|<\varepsilon_{\rm MR}.
\end{cases}
\label{eq:threshold-details}
\end{align}
The coarsest mean coefficient is always retained.
The thresholded state is then reconstructed as
\begin{align}
\bar{\boldsymbol S}^{\,\varepsilon}
=
\mathcal R_{\rm MR}
\left(
\bar S_0,\{d_{\ell,k}^{\,\varepsilon}\}
\right),
\label{eq:threshold-reconstruction}
\end{align}
where $\mathcal R_{\rm MR}$ denotes the inverse dyadic reconstruction.
The retained coefficient fraction is
\begin{align}
r_{\rm keep}
=
\frac{N_{\rm keep}}{N_{\rm total}},
\label{eq:retained-fraction}
\end{align}
with $N_{\rm keep}$ including the mean coefficient and all retained detail coefficients.
The thresholding error is measured relative to the unthresholded finite-volume state,
\begin{align}
e_j^\varepsilon
=
\bar S_j^{\,\varepsilon}
-
\bar S_j^{\,\mathrm{FV}},
\label{eq:threshold-error}
\end{align}
using the same RMSE, $L^1$, and $L^\infty$ metrics used later for the reference comparison.

The global mass change caused by thresholding is measured by
\begin{align}
\Delta M_\varepsilon
=
\left|
\sum_{j=1}^{N}
\left(
\bar S_j^{\,\varepsilon}
-
\bar S_j^{\,\mathrm{FV}}
\right)
\Delta x
\right|.
\label{eq:threshold-mass-defect}
\end{align}
Because the mean coefficient is retained, this quantity is expected to remain close to floating-point and reconstruction accuracy.
This diagnostic tests the conservation compatibility of the representation under controlled coefficient removal.
It does not replace the conservative finite-volume balance law.

Finally, the dyadic details are used to define a local adaptive indicator.
For front localization, the coarse levels are not the most useful because their support is broad.
We therefore form the indicator from the finest set of dyadic levels, denoted by $\mathcal L_{\rm fine}$:
\begin{align}
\eta_j
=
\sum_{\ell\in\mathcal L_{\rm fine}}
\sum_k
d_{\ell,k}^2\,
\mathbf 1_{j\in I_{\ell,k}},
\label{eq:local-adaptive-indicator}
\end{align}
where $I_{\ell,k}$ is the cell block associated with detail coefficient $d_{\ell,k}$.
The normalized indicator is
\begin{align}
\eta_j^{\rm norm}
=
\frac{\eta_j}{\max_{1\le m\le N}\eta_m}.
\label{eq:normalized-adaptive-indicator}
\end{align}
Cells with the largest values of $\eta_j$ are the locations where the saturation profile contains the strongest fine-scale activity.
In the numerical experiments, the top $5\%$ most active interior cells are marked as a prototype adaptive-refinement set.
This is only an indicator calculation on a fixed grid; it does not yet constitute an adaptive transport update.

The resulting workflow is therefore direct.
The saturation is advanced by the conservative finite-volume scheme.
The accepted cell-average state is projected into a bounded-interval multiwavelet representation and reconstructed back to cell averages.
The reconstructed state is used to measure FV--multiwavelet fidelity, while the dyadic detail hierarchy is used to quantify level activity, compression behavior, mass change under thresholding, and front localization.

\section{Numerical validation}
\label{sec:numerics}

\begin{figure}[t]
    \centering
    \includegraphics[width=\columnwidth]{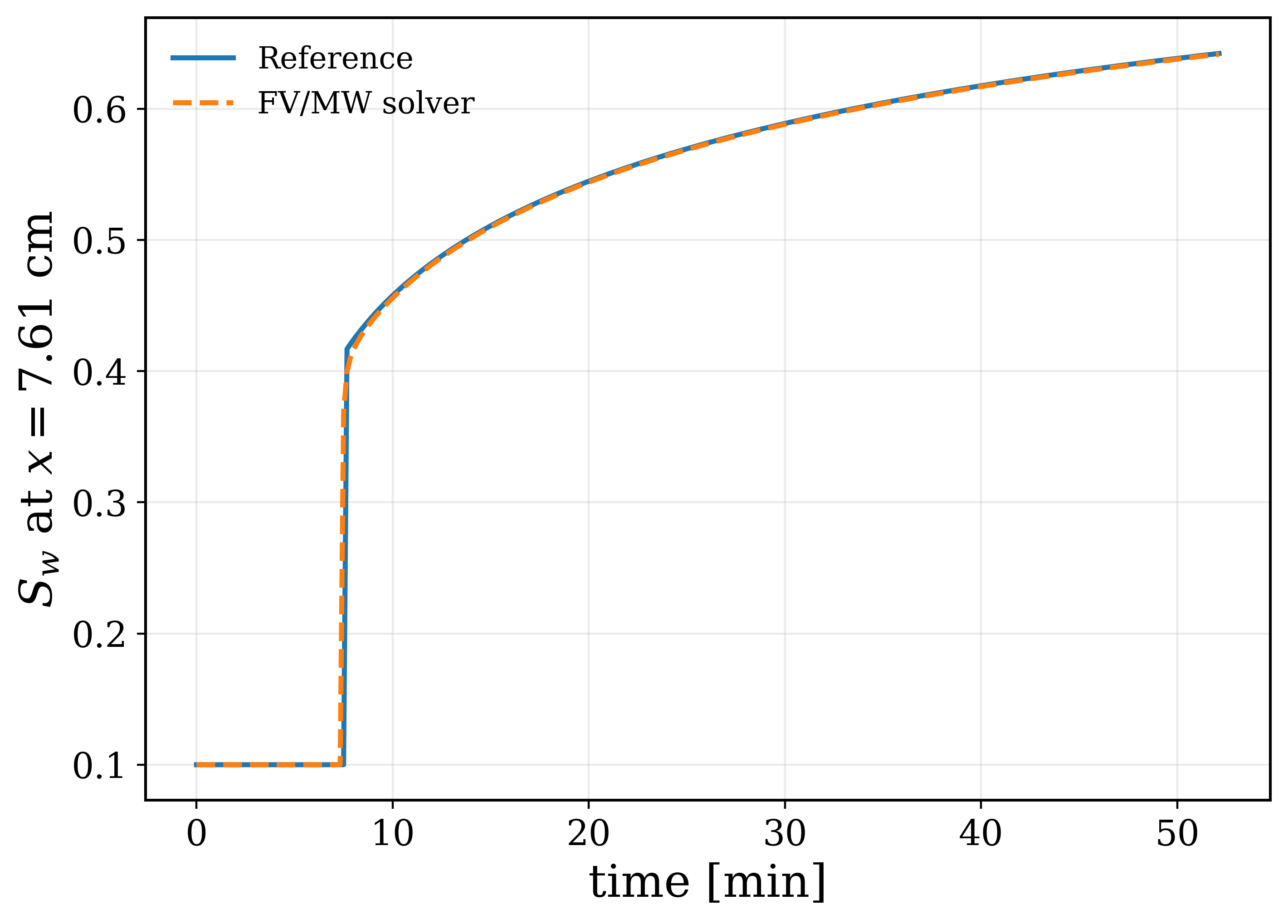}
\caption{Saturation history at the probe location $x=7.62\ \mathrm{cm}$.
The reference Buckley--Leverett solution is compared with the finite-volume
solution reconstructed through the bounded-interval multiwavelet representation.
The numerical profile reproduces the breakthrough jump and the post-front
increase in saturation, while the multiwavelet reconstruction does not alter
the local saturation history.}
    \label{fig:probe_saturation}
\end{figure}

\begin{figure}[t]
    \centering
    \includegraphics[width=\columnwidth]{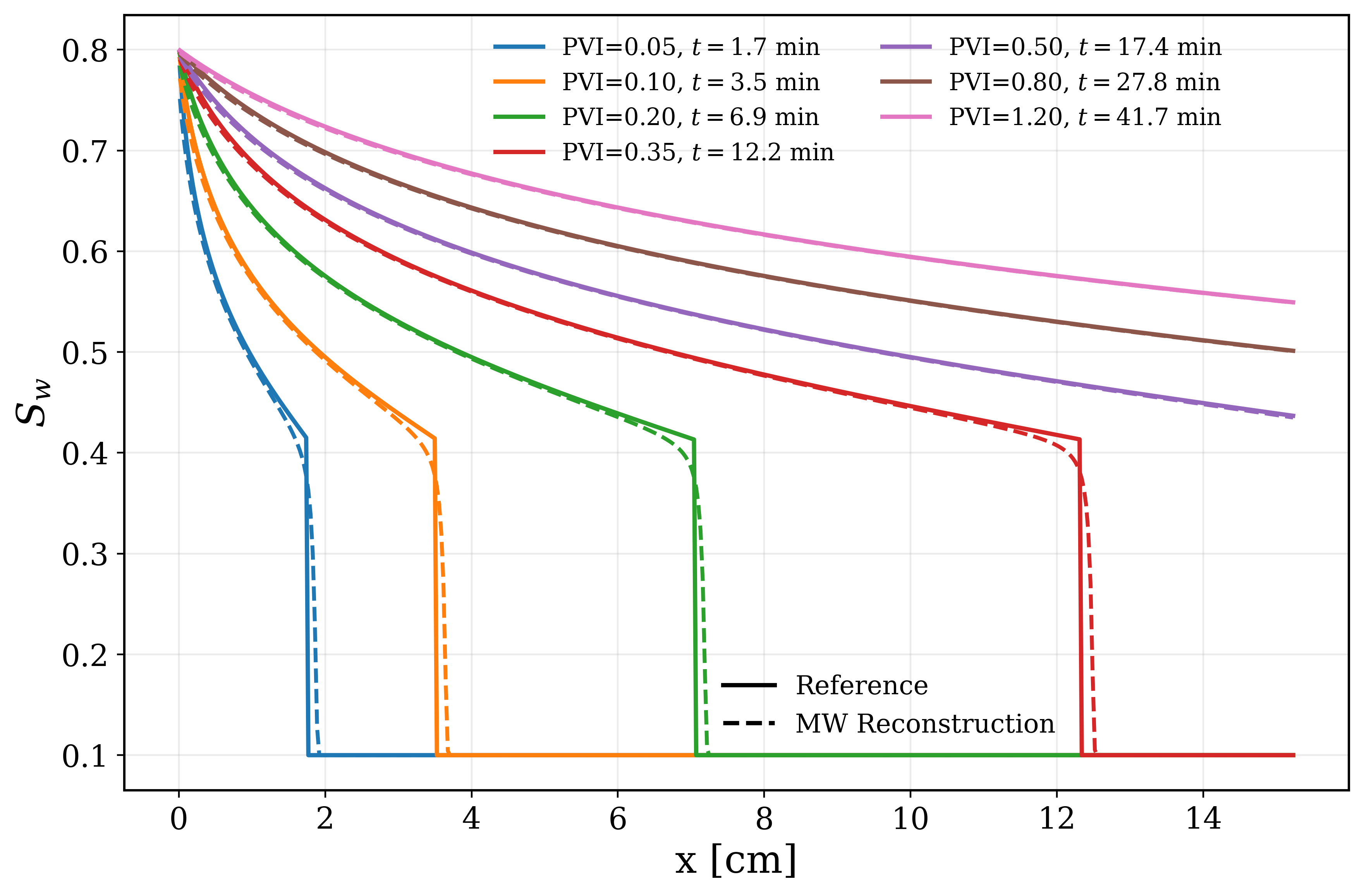}
    \caption{Spatial saturation profiles $S_w(x,t)$ at selected pore volumes injected (PVI). Solid lines denote the reference Buckley--Leverett profiles and dashed lines denote the bounded-interval MW-reconstructed numerical solution.}
    \label{fig:profiles_pvi}
\end{figure}

\begin{figure}[t]
    \centering
    \includegraphics[width=\columnwidth]{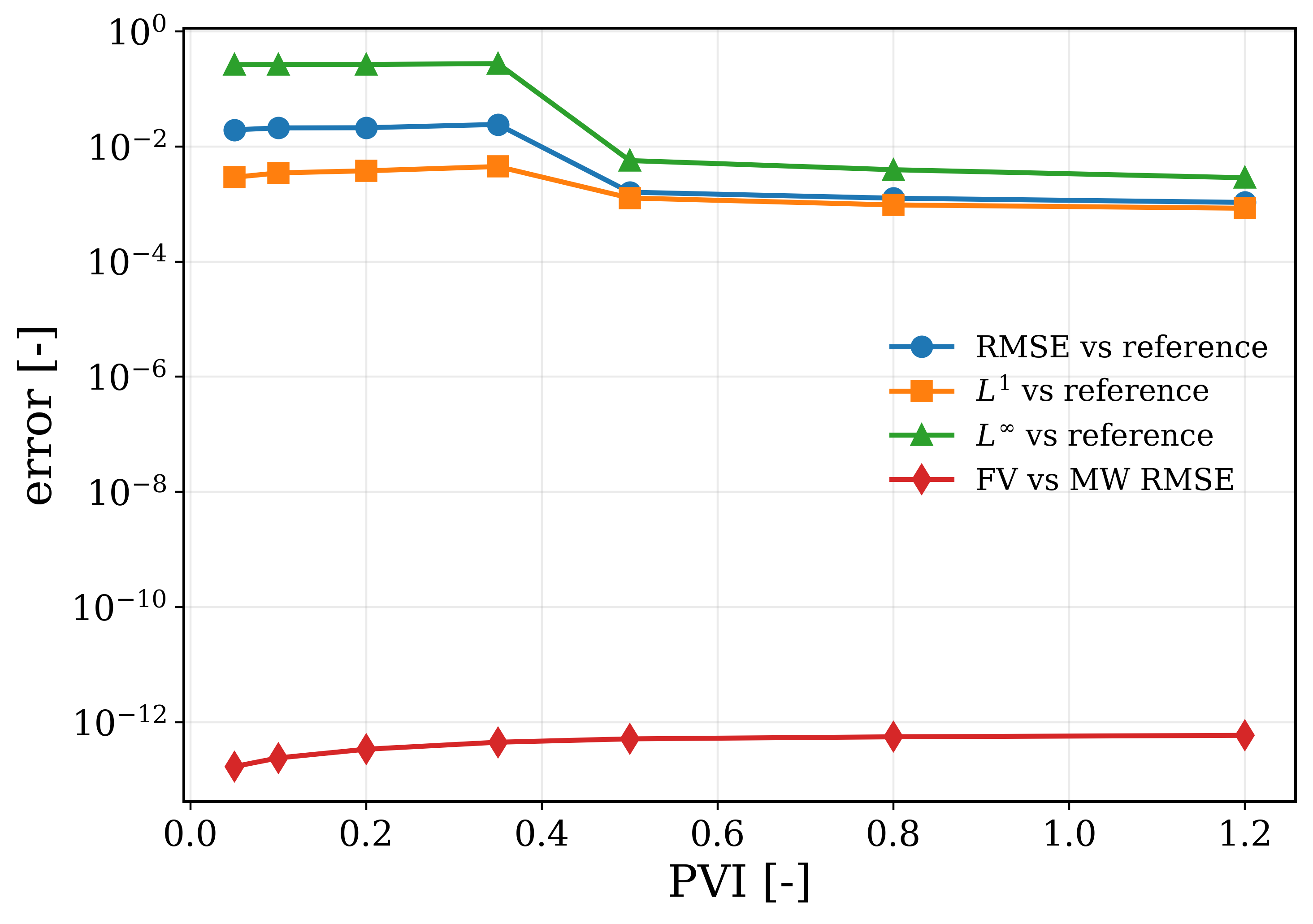}
    \caption{Snapshot-based error measures between the MW-reconstructed numerical state and the reference Buckley--Leverett profile as functions of pore volumes injected. Shown are the RMSE, the mean absolute error $L^1$, the maximum pointwise error $L^\infty$, and the internal consistency error between the finite-volume state and the MW reconstruction.}
    \label{fig:error_metrics}
\end{figure}

\begin{figure}[t]
    \centering
    \includegraphics[width=\columnwidth]{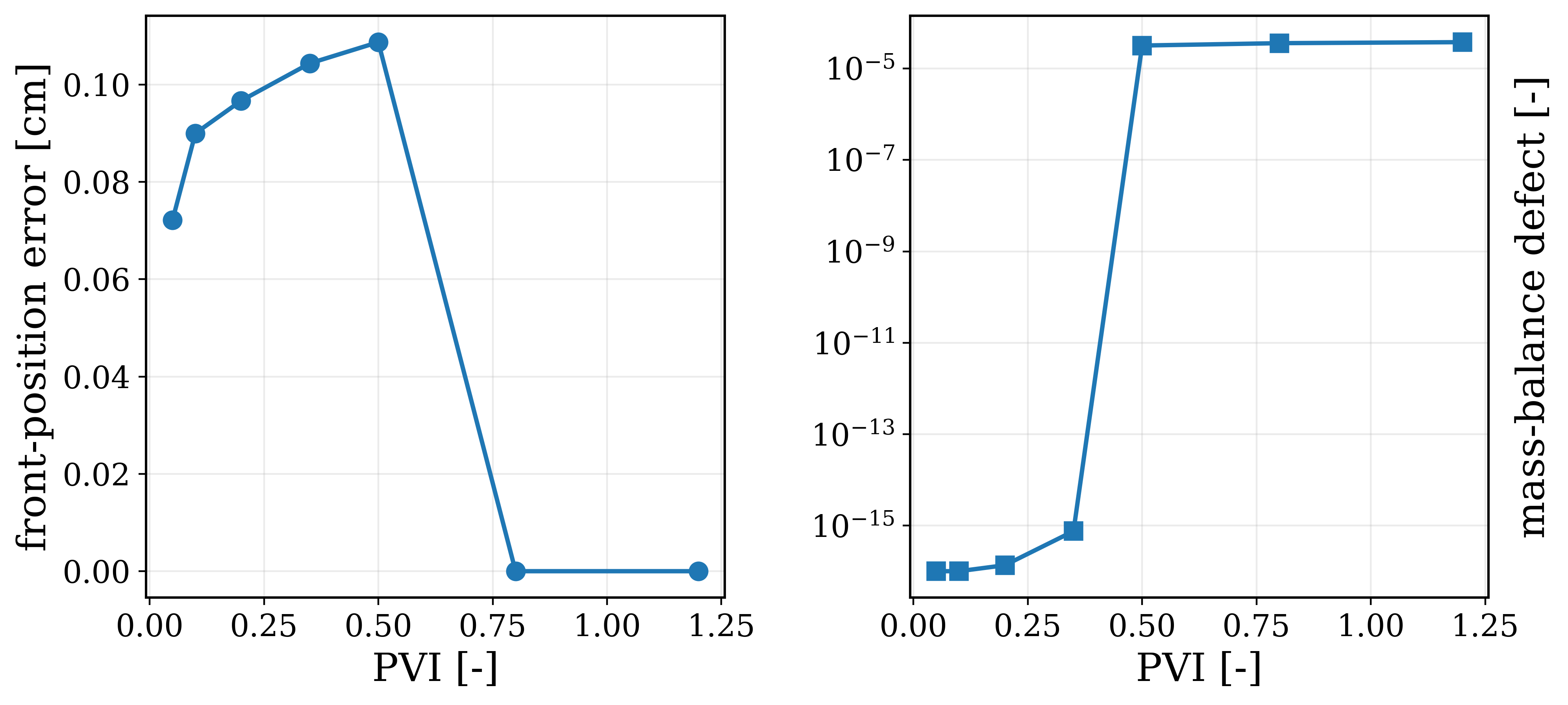}
    \caption{Additional transport diagnostics as functions of pore volumes injected. Left: absolute front-position error between the numerical and reference profiles, computed from a fixed saturation-threshold criterion. Right: mass-balance defect obtained by comparing the observed change in total water content with the time-integrated net boundary flux.}
    \label{fig:front_mass}
\end{figure}

\begin{figure}[t]
    \centering
    \includegraphics[width=\columnwidth]{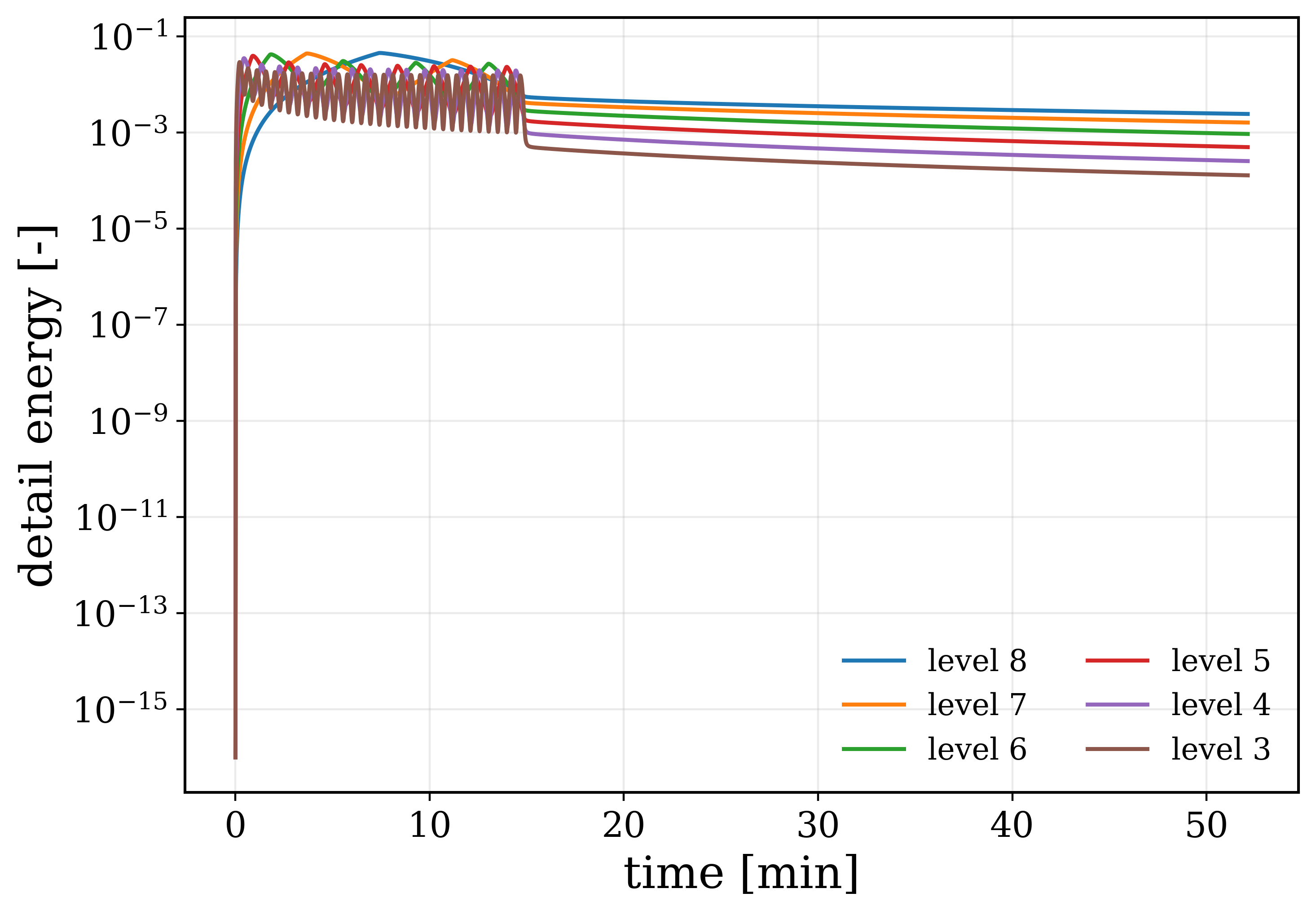}
    \caption{Temporal evolution of the most active dyadic detail energies $E_\ell(t)$ extracted from the reconstructed saturation state. Here, level $\ell$ denotes a dyadic multiresolution level obtained from repeated coarse--fine splitting of the state, and the plotted quantity is the corresponding detail energy.}
    \label{fig:detail_energies}
\end{figure}

\begin{figure*}[t]
    \centering
    \includegraphics[width=0.95\textwidth]{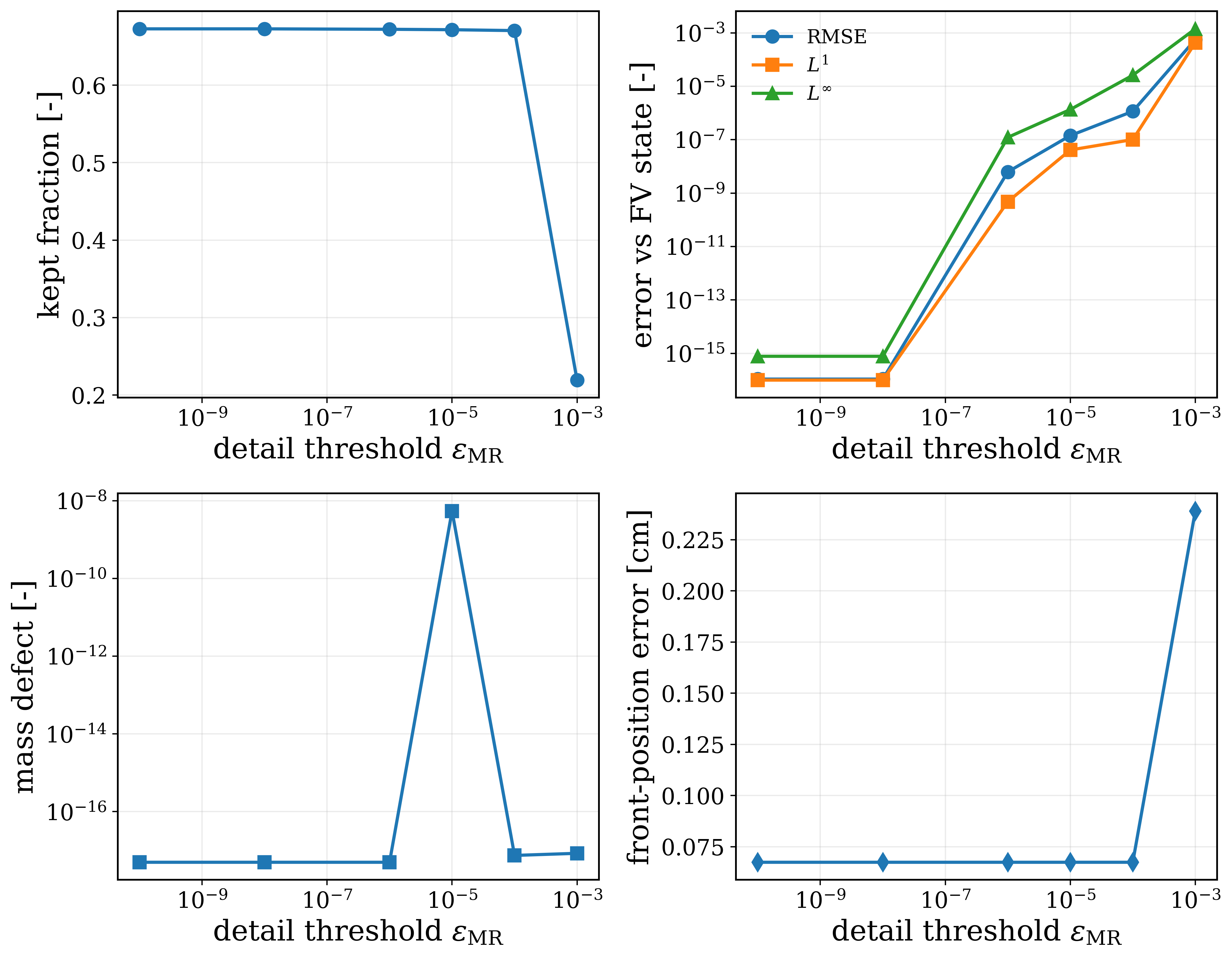}
\caption{Compression and mass-preservation behavior of the dyadic
multiresolution representation applied to accepted conservative finite-volume
states. Detail coefficients below $\varepsilon_{\rm MR}$ are discarded while
the global mean coefficient is retained. The retained coefficient fraction
decreases as the threshold is increased, whereas the reconstruction errors
remain small over a broad range of thresholds.
The mass defect remains close to roundoff, confirming that this thresholding experiment preserves the global mass content when the mean coefficient is retained.
The front-position
error increases only for the most aggressive threshold, where the shock
representation begins to degrade.}
    \label{fig:compression_conservation}
\end{figure*}

\begin{figure*}[t]
    \centering
    \includegraphics[width=0.95\textwidth]{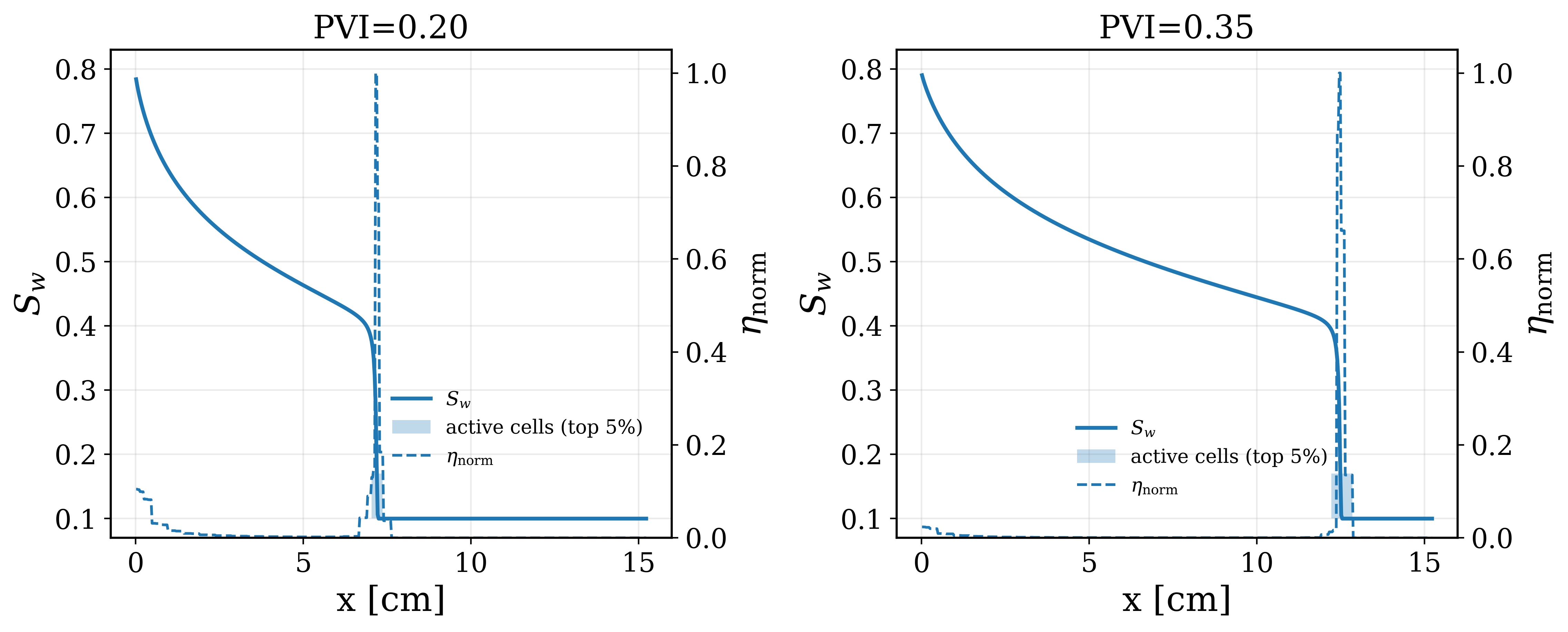}
    \caption{Fine-level multiresolution activity used as an adaptive indicator. The solid curve shows the finite-volume saturation profile and the dashed curve shows the normalized fine-level detail activity $\eta_{\rm norm}$. The shaded bands mark the top $5\%$ most active interior cells after excluding a $0.5~\mathrm{cm}$ boundary buffer from the marking step. The selected cells are localized around the moving Buckley--Leverett front, showing that the multiresolution representation supplies the spatial information needed for a future adaptive refinement strategy.}
    \label{fig:adaptive_indicator}
\end{figure*}

\begin{figure}[t]
    \centering
    \includegraphics[width=\columnwidth]{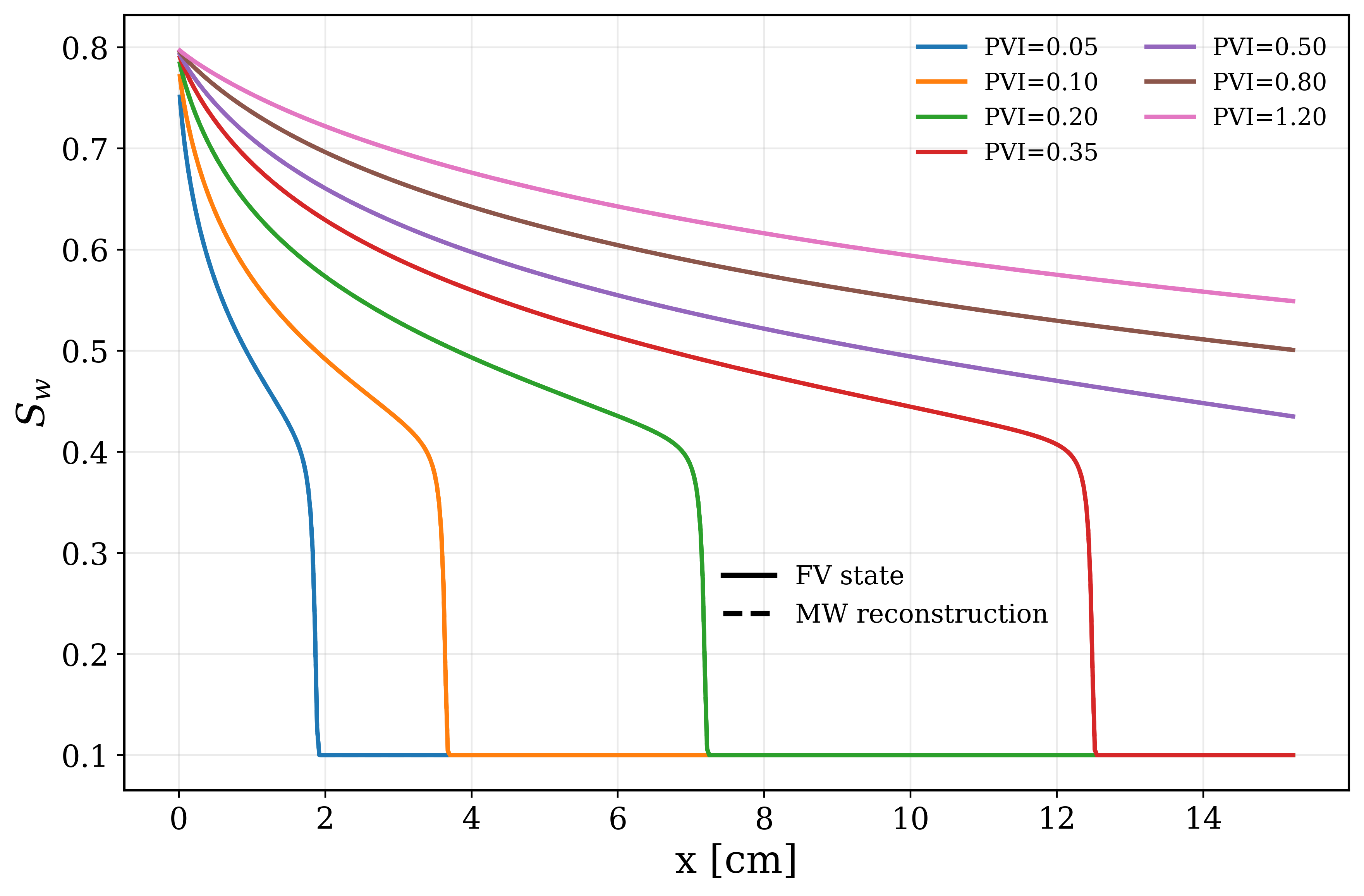}
    \caption{Direct comparison between the internal conservative finite-volume state and the bounded-interval MW-reconstructed state at selected pore volumes injected. The two profiles are practically indistinguishable over the full domain, with only extremely small differences near the steepest front.}
    \label{fig:fv_mw_internal}
\end{figure}

We assess the finite-volume/multiwavelet formulation on a deterministic
one-dimensional Berea benchmark. Reference profiles are generated with the
standard Buckley--Leverett fractional-flow construction implemented in
\texttt{pywaterflood}, using the same Corey-type constitutive parameters listed
in Table~\ref{tab:berea_setup}.\cite{Male2024pywaterflood}
The reference profiles are sampled on the computational grid and compared with
both the finite-volume state and the bounded-interval multiwavelet
reconstruction.

The validation addresses four points. First, the finite-volume transport
backbone should reproduce the reference saturation histories, spatial profiles,
front motion, and mass balance. Second, the multiwavelet reconstruction should
remain faithful to the accepted finite-volume state. Third, the dyadic detail
energies should expose the hierarchical activity generated by the moving front.
Fourth, coefficient thresholding and fine-level detail indicators should
quantify compression and front localization. Unless otherwise stated, all
results use the default benchmark parameters and solver settings in
Table~\ref{tab:berea_setup}.

\begin{table*}[t]
\caption{Physical, constitutive, and numerical parameters used in the default Berea benchmark calculations reported in this work, as extracted from the reference implementation.}
\label{tab:berea_setup}
\begin{ruledtabular}
\begin{tabular}{lll}
Quantity & Symbol & Value \\
\hline
\multicolumn{3}{l}{\textit{Physical and constitutive parameters}}\\
Core length & $L$ & $0.1524~\mathrm{m}$ ($15.24~\mathrm{cm}$) \\
Core diameter & $D$ & $0.0381~\mathrm{m}$ ($3.81~\mathrm{cm}$) \\
Cross-sectional area & $A=\pi D^2/4$ & $1.14009\times 10^{-3}~\mathrm{m^2}$ \\
Porosity & $\phi$ & $0.20$ \\
Connate water saturation & $S_{wc}$ & $0.10$ \\
Residual oil saturation & $S_{or}$ & $0.20$ \\
Initial water saturation & $S_{w,\mathrm{init}}$ & $0.10$ \\
Injected water saturation & $S_{w,\mathrm{inj}}$ & $0.80$ \\
Water viscosity & $\mu_w$ & $1.0\times 10^{-3}~\mathrm{Pa\,s}$ \\
Oil viscosity & $\mu_o$ & $4.0\times 10^{-3}~\mathrm{Pa\,s}$ \\
Corey exponent (water) & $n_w$ & $2.0$ \\
Corey exponent (oil) & $n_o$ & $2.0$ \\
Endpoint water relative permeability & $k_{rw}^0$ & $1.0$ \\
Endpoint oil relative permeability & $k_{ro}^0$ & $1.0$ \\
Injection rate & $q$ & $1.0~\mathrm{mL/min}$ \\
Darcy velocity & $v=q/A$ & $1.26306~\mathrm{m/day}$ \\
\hline
\multicolumn{3}{l}{\textit{Numerical and multiwavelet parameters}}\\
Number of finite-volume cells & $N$ & $512$ \\
Grid spacing & $\Delta x=L/N$ & $2.97656\times10^{-4}~\mathrm{m}$ \\
Numerical flux & --- & Godunov \\
Time integrator & --- & SSPRK2 \\
CFL number & $\mathrm{CFL}$ & $0.85$ \\
Final simulated pore volumes injected & --- & $1.50$ \\
Final simulation time & $t_{\mathrm{end}}$ & $0.03620~\mathrm{day}$ ($52.125~\mathrm{min}$) \\
Snapshot PVI values & --- & $0.05,\ 0.10,\ 0.20,\ 0.35,\ 0.50,\ 0.80,\ 1.20$ \\
Probe location & $x_p$ & $L/2 = 0.0762~\mathrm{m}$ ($7.62~\mathrm{cm}$) \\
Bounded-interval MW order & --- & $8$ \\
MW projection precision & --- & $10^{-7}$ \\
MW reconstruction quadrature & --- & 8-point Gauss--Legendre \\
Front-location threshold & --- & $0.5$ \\
\end{tabular}
\end{ruledtabular}
\end{table*}

The benchmark uses a one-dimensional core with fixed porosity, Corey-type
relative permeabilities, constant total Darcy velocity, and non-capillary
Buckley--Leverett transport. The spatial domain is discretized by a uniform
finite-volume mesh with $N=512$ cells, which is also convenient for the dyadic
multiresolution analysis. The transport step uses SSPRK2 time integration and
the Godunov numerical flux. The multiwavelet layer is applied only to accepted
finite-volume states: the state is projected into the bounded-interval
representation and reconstructed back to cell averages for comparison,
compression, mass-defect diagnostics, and adaptive-indicator analysis.

A first validation quantity is the saturation history at a fixed probe location.
Let $x_p$ denote the probe position, chosen here near the midpoint of the core.
We compare the numerical prediction
\begin{align}
t \mapsto S_w(x_p,t)
\end{align}
against the corresponding reference Buckley--Leverett value.
The result is shown in Figure~\ref{fig:probe_saturation}.
This figure is particularly important because it condenses the transport behavior into a single physically interpretable observable.
Before the front arrives, the probe saturation must remain near the initial saturation.
At breakthrough, the probe must register the correct sharp increase.
After breakthrough, the solver must continue to reproduce the gradual rise associated with the post-front fractional-flow dynamics.
The agreement in Figure~\ref{fig:probe_saturation} is excellent over the full time interval.
The breakthrough jump occurs at essentially the same time in the numerical and reference solutions, and the subsequent evolution remains almost indistinguishable.
This confirms that the conservative finite-volume backbone reproduces the correct one-dimensional front propagation and that the added multiwavelet representation does not corrupt the physically relevant local saturation history.

We next turn to the full spatial saturation profiles at selected pore volumes injected.
For a prescribed sequence of PVI values, we extract the numerical profiles
\begin{align}
x \mapsto S_w(x,t_{\mathrm{PVI}})
\end{align}
and compare them with the corresponding reference Buckley--Leverett curves.
The results are collected in Figure~\ref{fig:profiles_pvi}.
Several features deserve attention.
First, the front positions are reproduced correctly over the entire sequence of snapshots, from early-time entry of the front into the core to later-time propagation across a large fraction of the domain.
Second, the numerical profiles preserve the expected Buckley--Leverett structure: a monotone trailing branch behind the shock and a sharp transition to the initial saturation ahead of the front.
Third, the bounded-interval MW reconstruction tracks the reference profile extremely closely.
The visible discrepancy is confined mainly to the steepest shock region, where any discrete representation of a nearly discontinuous profile is naturally most sensitive.
Even there, the mismatch remains small.
Thus Figure~\ref{fig:profiles_pvi} already shows at the level of the profiles themselves that the present formulation captures the correct deterministic Buckley--Leverett dynamics while keeping the state in a bounded-interval multiresolution representation.

To move beyond visual comparison, we evaluate snapshot-based errors between the MW-reconstructed numerical state and the reference profile.
Defining the cellwise snapshot error by
\begin{align}
e_j(t_n)=\bar S_j^{\,\mathrm{MW}}(t_n)-\bar S_j^{\,\mathrm{ref}}(t_n),
\end{align}
we compute the root-mean-square error, the discrete mean absolute error, and the discrete maximum error as
\begin{align}
\mathrm{RMSE}(t_n)
&=
\left(
\frac{1}{N}
\sum_{j=1}^{N}
\left(
\bar S_j^{\,\mathrm{MW}}(t_n)-\bar S_j^{\,\mathrm{ref}}(t_n)
\right)^2
\right)^{1/2},
\\
L^1(t_n)
&=
\frac{1}{N}
\sum_{j=1}^{N}
\left|
\bar S_j^{\,\mathrm{MW}}(t_n)-\bar S_j^{\,\mathrm{ref}}(t_n)
\right|,
\\
L^\infty(t_n)
&=
\max_{1\le j\le N}
\left|
\bar S_j^{\,\mathrm{MW}}(t_n)-\bar S_j^{\,\mathrm{ref}}(t_n)
\right|.
\label{eq:error-metrics}
\end{align}

In addition, to measure the fidelity of the representation layer itself, we compare the internal conservative finite-volume state and the reconstructed multiwavelet state through
\begin{align}
\mathrm{RMSE}_{\mathrm{FV\text{-}MW}}(t_n)
=
\left(
\frac{1}{N}
\sum_{j=1}^{N}
\left(
\bar S_j^{\,\mathrm{FV}}(t_n)-\bar S_j^{\,\mathrm{MW}}(t_n)
\right)^2
\right)^{1/2}.
\end{align}
These quantities are reported in Figure~\ref{fig:error_metrics}.

The interpretation of Figure~\ref{fig:error_metrics} is important.
The FV--MW reconstruction error is essentially at machine precision over the entire set of snapshots, which demonstrates that the bounded-interval projection and reconstruction act as a numerically faithful representation of the conservative state.
In other words, once the transport step has produced the finite-volume state, the multiwavelet layer does not distort it in any measurable way.
By contrast, the errors relative to the external Buckley--Leverett reference remain finite, as expected, because the transport solution is still computed on a discrete finite-volume grid.
The $L^1$ and RMSE values remain small throughout the simulation, indicating good global profile agreement.
The $L^\infty$ error is largest at early and intermediate PVI because this metric is dominated by the steepest local discrepancy, which in a shock-dominated problem is naturally concentrated near the front.
A small positional mismatch in the shock region can therefore produce a noticeably larger $L^\infty$ value even when the full profile is otherwise well matched.
This is precisely why the simultaneous presentation of RMSE, $L^1$, and $L^\infty$ is useful:
the first two reflect overall profile agreement, while the last isolates the worst pointwise error near the steep front.
Taken together, the three measures show that the numerical solution remains globally accurate, while the multiwavelet representation itself is essentially exact relative to the internal conservative state.
This interpretation is also consistent with the direct comparison between the internal finite-volume state and the reconstructed multiwavelet state, shown in Figure~\ref{fig:fv_mw_internal}, where the two are practically indistinguishable over the full domain except for extremely small differences near the steepest front.

A further physically meaningful diagnostic is the front-position error.
At each snapshot we estimate the front location in both the numerical and reference profiles through a fixed saturation threshold and define the absolute difference between the two positions.
At the same time, because the transport update is conservative, we monitor the mass-balance defect obtained by comparing the observed change in total water content with the time-integrated net boundary flux.
These two diagnostics are shown in Figure~\ref{fig:front_mass}.
The front-position error remains small throughout the simulation and is at worst only a small fraction of the core length.
Its early-time increase reflects the greater sensitivity of front localization when the shock is still entering and steepening inside the domain.
At later PVI values the front-position discrepancy drops further, and for the latest snapshots it becomes essentially negligible on the scale of the plotted domain.
The mass-balance defect remains extremely small, close to machine precision at early times and still negligible at later times.
The present agreement with the Buckley--Leverett reference is not obtained by sacrificing conservative structure; on the contrary, it is achieved precisely because the transport backbone retains the local conservative flux form that enforces the correct shock dynamics.

Beyond direct validation against the reference solution, the present formulation also provides access to the multiresolution structure of the evolving state.
To quantify that structure, we compute the dyadic detail energies introduced in Section~\ref{sec:optionA},
\begin{align}
E_\ell(t)=\sum_k d_{\ell,k}(t)^2,
\end{align}
where $d_{\ell,k}(t)$ denotes the dyadic detail coefficient at level $\ell$ extracted from the reconstructed state.
Figure~\ref{fig:detail_energies} displays the temporal evolution of the most active levels.
In this figure, the label ``level'' refers to the dyadic multiresolution level of the diagnostic coarse--fine splitting and not to a physical energy level.
These curves should not be interpreted as validation errors; rather, they are structural diagnostics of how the solution populates the hierarchy.
At very early times, as the front begins to form and enter the computational domain, several levels exhibit a rapid growth in detail energy, indicating that the solution develops sharp local structure.
During the interval in which the shock traverses the core most actively, the energies remain elevated across a band of levels, reflecting the presence of strong multiscale content.
At later times, as the transported profile becomes more spatially extended and relatively smoother away from the sharp transition, the detail energies decay gradually.
The different levels decay at different rates, which is consistent with the interpretation that fine- and intermediate-scale contributions respond differently to the motion of the front.
This figure shows that the moving front generates a measurable hierarchy of
dyadic activity. The next two diagnostics use this hierarchy more directly for
coefficient thresholding and front localization.

We next perform a thresholding experiment on accepted finite-volume states.
This test does not alter the transport update. Its purpose is to determine how
much dyadic detail content can be removed while preserves the global mass content. Detail coefficients below
$\varepsilon_{\rm MR}$ are discarded, the global mean coefficient is retained,
and the thresholded state is reconstructed back to finite-volume cell averages.

Figure~\ref{fig:compression_conservation} reports the retained coefficient
fraction, reconstruction errors, mass defect, and front-position error. For
small and moderate thresholds, the RMSE and $L^1$ errors remain small while the
retained coefficient fraction decreases, showing that the multiresolution
representation contains compressible detail content. The $L^\infty$ error is
more sensitive because it is controlled by the sharp front. At the most
aggressive threshold, $\varepsilon_{\rm MR}=10^{-3}$, the errors and
front-position discrepancy increase, indicating that dynamically relevant
front information has begun to be removed.

The mass defect remains close to roundoff over the threshold range. This follows
from retaining the coarsest mean coefficient. The result should not be
interpreted as a statement that arbitrary coefficient manipulation is locally
conservative; local conservation is guaranteed by the finite-volume flux update.
Rather, Figure~\ref{fig:compression_conservation} shows that this particular
thresholding diagnostic preserves the global mass content while
quantifying the cost of coefficient removal.

Figure~\ref{fig:adaptive_indicator} shows the fine-level detail activity used
as a local front indicator. The normalized indicator $\eta_{\rm norm}$ peaks at
the moving Buckley--Leverett front, and the marked cells, defined as the top
$5\%$ most active interior cells after excluding a narrow boundary buffer from
the visualization, are localized around the front at both reported PVI values.
This is the expected behavior of a multiresolution indicator: smooth regions
carry weak fine-level detail activity, whereas the steep transition generates a
localized signal. The calculation is not a fully adaptive transport simulation,
because the solution is still advanced on a fixed finite-volume grid. It
nevertheless demonstrates that the bounded-domain multiresolution layer
provides the spatial activity information needed for a later adaptive extension.

Taken together, Figures~\ref{fig:probe_saturation}--\ref{fig:adaptive_indicator}
support the intended interpretation of the method. The finite-volume backbone
reproduces the reference Buckley--Leverett dynamics, including breakthrough
behavior, front location, and mass balance. The bounded-interval multiwavelet
layer reconstructs the accepted finite-volume state with essentially exact
fidelity and exposes the hierarchical organization of the transported profile.
The thresholding and adaptive-indicator diagnostics show that the hierarchy also
provides controlled compression, global mass-preservation diagnostics, and
localized front-activity information.
These results establish the present formulation as a conservative
finite-volume Buckley--Leverett solver equipped with a bounded-domain
multiresolution analysis layer.

\section{Conclusions}
\label{sec:conclusions}

We have developed and validated a conservative finite-volume Buckley--Leverett solver equipped with a bounded-interval multiwavelet state-analysis layer.
The formulation separates transport from representation.
The nonlinear hyperbolic saturation equation is advanced by a finite-volume method in conservative flux-difference form, while the accepted saturation state is projected into a bounded-domain multiwavelet hierarchy and reconstructed back to cell averages for analysis.

The Berea benchmark results show that this separation is effective.
The conservative finite-volume backbone reproduces the reference Buckley--Leverett saturation histories, spatial profiles, front positions, and mass balance.
The bounded-interval multiwavelet reconstruction tracks the accepted finite-volume state with essentially exact fidelity, confirming that the representation layer does not distort the transported saturation profile when used diagnostically.

The multiresolution diagnostics clarify the practical role of the representation.
The detail energies expose the hierarchical activity generated by the moving front.
The thresholding experiments show that a substantial fraction of dyadic detail coefficients can be removed while retaining small reconstruction errors and negligible global mass defect when the mean coefficient is preserved.
The fine-level detail activity also localizes the moving Buckley--Leverett front, providing a prototype adaptive indicator on the fixed finite-volume grid.
Thus, the multiwavelet layer supplies controlled compression, mass-preservation diagnostics, and front-localization information beyond direct visualization of the finite-volume solution.

The present method should be interpreted as a conservative first stage toward transport-active multiresolution algorithms, not as a fully native multiwavelet transport solver.
Because the transport update is performed by the finite-volume operator, the formal convergence behavior remains that of the underlying finite-volume discretization.
Future work will build on the diagnostics established here by developing conservative multilevel transfer, adaptive thresholding with local cell-average correction, and interface-flux mechanisms operating across a dynamic multiwavelet hierarchy.

\section*{Data availability}
The data supporting the findings of this study are available from the corresponding author upon reasonable request.

\section*{Code availability}
The code used in this work is publicly available at \url{https://github.com/Christian48596/fv_mw_buckley_leverett}. The repository contains the Buckley--Leverett finite-volume/bounded-interval multiwavelet solver, the default Berea benchmark settings, and documentation for modifying the rock and fluid parameters. The results reported in this study can be reproduced from the code and input settings provided in that repository.

\section*{Funding}
This work was supported by the Deanship of Research (DOR) at King Fahd University of Petroleum \& Minerals (KFUPM) under Project No.~EC251017.

\section*{Acknowledgements}
Ch.~T. thanks Evgueni Dinvay (UiT) for useful discussions.

\section*{Declarations}

\subsection*{Competing interests}
The author declares no competing interests.

\subsection*{Ethics approval}
Not applicable.

\subsection*{Consent to participate}
Not applicable.

\subsection*{Consent for publication}
The author consents to publication of this manuscript.

\bibliography{main}

\end{document}